# On the Tight Chromatic Bounds for a Class of Graphs without Three Induced Subgraphs

Medha Dhurandhar

**Abstract:** Here we prove that a graph without some three induced subgraphs has chromatic number at the most equal to its maximum clique size plus one.

**Introduction:**

It has been an eminent unsolved problem in graph theory to determine the chromatic number of a given graph. Failing in the efforts to determine this, attempts have been made to find good bounds for the chromatic number of a graph. Vizing [4] proved that if a graph does not induce some nine subgraphs, then $\omega(G) \leq \chi(G) \leq \omega(G)+1$ where $\omega(G)$ is the size of maximum clique in G and $\chi(G)$ is the chromatic number of G. Later Choudum [1] and Javdekar [2] improved this result by dropping five and six of these nine subgraphs from the hypothesis, respectively. Finally Kierstead [3] showed that $\omega(G) \leq \chi(G) \leq \omega(G)+1$ for a $\{K_{1,3},(K_5-e)\}$-free graph. Furthermore, Dhurandhar [5] proved that $\omega(G) \leq \chi(G) \leq \omega(G)+1$ for a $\{K_{1,3}, (K_2 \cup K_1)+K_2\}$-free graph. In this paper we prove that if G is $\{K_{1,3}, H_1, H_2\}$-free, then $\chi(G) \leq \omega(G)+1$ where

$H_1 =$ 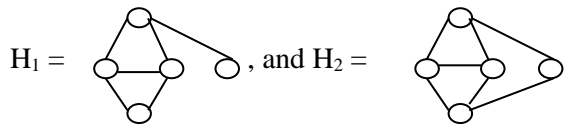 , and $H_2 =$

**Notation:** For a graph G, V(G), E(G), $\Delta(G)$, $\omega(G)$, $\chi(G)$ denote the vertex set, edge set, maximum degree, size of a maximum clique, chromatic number respectively. For $u \in V(G)$, $N(u) = \{v \in V(G)$ / $uv \in E(G)\}$, and $\overline{N(u)} = N(u) \cup (u)$. If $S \subseteq V(G)$, then $<S>$ denotes the subgraph of G induced by S. Also for $u \in V(G)$, deg u and $<Q_u>$ denote the degree of u in G and a maximum clique in N(u). If C is some coloring of G and if a vertex u of G is colored m in C, then u is called a m-vertex. All graphs considered henceforth are simple.

First we prove a Lemma, which will be used later.

**Lemma 1:** Let G be $\{K_{1,3}, H_1, H_2\}$-free and $v \in V(G)$. Then either

1. $<N(v)-Q_v>$ is complete or
2. $\{N(v)-Q_v\} = \{w, z\}$, $zw \notin E(G)$. Further $\exists$ w', z' $\in Q_v$ s.t. ww', zz' $\notin E(G)$ and $zq \in E(G)$ $\forall$ $q \in Q_v-z'$, $wq \in E(G)$ $\forall q \in Q_v-w'$ or
3. $\{N(v)-Q_v\} = \{w, z, x\}$ with xz, xw $\in E(G)$ and wz $\notin E(G)$. Further $\exists$ w', z' $\in Q_v$ s.t. ww', zz', xw', xz' $\notin E(G)$ and wq $\in E(G)$ $\forall q \in Q_v-w'$, zq $\in E(G)$ $\forall q \in Q_v-z'$, xq $\in E(G)$ $\forall q \in Q_v-\{w', z'\}$.

Proof: Let $\exists$ $v \in V(G)$ with w, z $\in N(v)-Q_v$ s.t. wz $\notin E(G)$. Then $\exists$ w', z' $\in Q_v$ s.t. ww', zz' $\notin E(G)$. As G is $K_{1,3}$-free, wz', zw' $\in E(G)$ and $\forall$ q $\in Q_v-\{w', z'\}$, qw $\in E(G)$ or qz $\in E(G)$. W.l.g. let qw $\in$ E(G). Then qz $\in$ E(G) (else $<q, z', w, w', z> = H_1$). This proves 2, in case $|N(v)-Q_v| = 2$. Let $|N(v)-Q_v|$

> 2 and {N(v)-$Q_v$} ⊇ {w, z, x} with wz ∉ E(G). As before ∃ w', z' ∈ $Q_v$ s.t. ww', zz' ∉ E(G) and wq ∈ E(G) ∀ q ∈ $Q_v$-w', zq ∈ E(G) ∀ q ∈ $Q_v$-z'. W.l.g. let xw ∈ E(G).

**Claim:** xz ∈ E(G)

If not, then xz' ∈ E(G) and as ω<N(v)> ≥ 3, ∃ q ∈ $Q_v$-{w', z'}. Now qz ∈ E(G) iff qw ∈ E(G) (else <q, w', z, z', w> = $H_1$). Hence as G is $K_{1,3}$-free, qw, qz ∈ E(G). Again xq ∈ E(G) (else <w, z', q, x, z> = $H_1$). But then <$Q_v$-w'+{w, x}> is a bigger clique in N(V) than $Q_v$, a contradiction.

This proves the Claim.

Again xq ∈ E(G) ∀ q ∈ $Q_v$-{w', z'} (else xw' ∈ E(G) and <z, w', q, w, x> = $H_2$) and hence xw', xz' ∉ E(G). Now if |N(v)-$Q_v$| > 3, then let y ∈ N(v)-$Q_v$. As before yq ∈ E(G) ∀ q ∈ $Q_v$-{w', z'} and yw', yz' ∉ E(G). But then if xy ∈ E(G) <$Q_v$-{z', w'}+{w, x, y}> is a bigger clique in N(V) than $Q_v$ and if xy ∉ E(G) <v, w', x, y> = $K_{1,3}$, a contradiction. Hence |N(v)-$Q_v$| = 3 and this proves 3.

This completes the proof of **Lemma 1**.

For completeness we simply state a result from [6], which will be used in the main result.

**Lemma 2:** If G is $K_{1,3}$-free, then each component of the subgraph of G induced by two colour classes is either a path or a cycle.

**Theorem:** If G is a {$K_{1,3}$, $H_1$, $H_2$}-free, then χ(G) ≤ ω(G)+1.

Proof: Let G be a smallest {$K_{1,3}$, $H_1$, $H_2$}-free graph with χ(G) > ω(G)+1. Now by minimality ∀ u ∈ V(G), χ(G-u) ≤ ω(G-u)+1. Thus ω(G)+1 < χ(G) ≤ χ(G-u)+1 ≤ ω(G-u)+2 ≤ ω(G)+2. Hence ∀ u ∈ V(G), χ(G-u) = ω(G)+1 and hence χ(G) = ω(G)+2.

**Case 1:** ∃ v ∈ V(G) s.t. {N(v)-$Q_v$} = {w, z} with zw ∉ E(G).

Let w', z' ∈ $Q_v$ be s.t. ww', zz' ∉ E(G). Then by **Lemma 1,** zq ∈ E(G) ∀ q ∈ $Q_v$-z', wq ∈ E(G) ∀ q ∈ $Q_v$-w'. Let C = $\bigcup_{1}^{\omega+2} i$ be a (ω(G)+2)-coloring of G in which v receives the unique color ω+2 and vertices in $Q_v$ receive colors 1, ..., |$Q_v$|. W.l.g. let w be the ω-vertex and z be the (ω+1)-vertex of v.

Then by **Lemma 2**, ∃ a w-z path P s.t. vertices on P are alternately colored ω and ω+1 (else w can be colored by ω+1 and v by ω). Let x be the ω+1-vertex adjacent to w on P. Then xz' ∈ E(G) (else <z', v, w', w, x> = $H_1$ or $H_2$). But then <z', w, x, v, z> = $H_1$, a contradiction.

**Case 2:** ∃ v ∈ V(G) s.t. {N(v)-$Q_v$} = {w, z, x} with xz, xw ∈ E(G) and wz ∉ E(G).

Let w', z' ∈ $Q_v$ be s.t. ww', zz' ∉ E(G). Then by **Lemma 1**, xw', xz' ∉ E(G), wq ∈ E(G) ∀ q ∈ $Q_v$-w', zq ∈ E(G) ∀ q ∈ $Q_v$-z', and xq ∈ E(G) ∀ q ∈ $Q_v$-{w', z'}. Let C = $\bigcup_{1}^{\omega+2} i$ be a (ω(G)+2)-coloring of

G in which v receives the unique color $\omega+2$ and vertices in $Q_v$ receive colors 1, ..., $|Q_v|$. Clearly colors $\omega$ and $\omega+1$ are used in $N(v)-Q_v$ (else color v by the missing color).

**Case 2.1:** x has color say $\omega$.

Again $N(v)-Q_v$ has a $\omega+1$-vertex (else color v by $\omega+1$). Thus either w' or z' has a color not used in $N(v)-Q_v$. W.l.g. let z' have a color not used in $N(v)-Q_v$. As z'x $\notin$ E(G), z' has a $\omega$-vertex say y outside N(v). Then yw $\in$ E(G) (else <v, w, x, z', y> = $H_1$). But then <z', w, y, v, z> = $H_1$ or $H_2$, a contradiction.

**Case 2.2:** x has color other than $\omega$ and $\omega+1$.

Then w.l.g. let w, z be colored by $\omega$ and $\omega+1$ resply. As before either w' or z' has a color not used in $N(v)-Q_v$. W.l.g. let z' have a color not used in $N(v)-Q_v$. As z'z $\notin$ E(G), z' has a $\omega+1$-vertex say y outside N(v). Then yw $\in$ E(G) (else <v, w, x, z', y> = $H_1$ or $H_2$). But then <z', w, y, v, z> = $H_1$, a contradiction.

**Case 3:** $\forall$ v $\in$ V(G), <N(v)-$Q_v$> is complete

Let Q be a maximum clique in G and v $\in$ Q. Let $C = \bigcup_1^{\omega+2} i$ be a ($\omega(G)+2$)-coloring of G in which v receives the unique color $\omega+2$ and vertices in Q receive colors 1, ..., $\omega-1$. Label v as $v_0$. Clearly colors $\omega$ and $\omega+1$ are used in N(v)-Q (else color v by the missing color). Hence as $|N(v)-Q| \leq |Q|$, $\exists$ a color say 1 not used in N(v)-Q. Label the vertex in Q, which has color 1 as $v_1$ ($\exists$ such a vertex, else color v by 1). Again $\exists$ a color say 2 not used in $N(v_1)$-Q and a 2-vertex in Q (else color $v_1$ by 2 and v by 1). Label this 2-vertex in Q as $v_2$ and so on. Label v as $v_0$. Let $v_0, v_1,..., v_k$ be a maximal sequence of vertices in Q s.t. $v_i$ has color i and color i+1 is not used in $N(v_i)$-Q, $1 \leq i \leq k$. By maximality of the sequence k+1 = t+1 for some $0 \leq t \leq k-2$. Consider a component P containing $\omega$-vertex of $v_t$ s.t. vertices in P are colored t+1 or $\omega$. As $v_t$ has a unique $\omega$-vertex, by **Lemma 2**, P is a path. If $v_{t+1} \notin$ P, then alter colors of vertices in P, color $v_t$ by $\omega$ and $v_j$ by j+1 for $0 \leq j \leq t-1$, a contradiction. Hence $v_{t+1} \in$ P. Similarly if R is a component containing $\omega$-vertex of $v_k$ s.t. vertices in R are colored t+1 or $\omega$, then $v_{t+1} \in$ R. Thus P=R and as G is $K_{1,3}$-free at least two of $v_t, v_{t+1}, v_k$ have the same $\omega$-vertex. Similarly at least two of $v_t, v_{t+1}, v_k$ have the same $\omega+1$-vertex.

**Case 3.1:** At the most two of $v_t, v_{t+1}, v_k$ have the same $\omega$-vertex.
Let {m, n, p} = {t, t+1, k} and $v_m, v_n$ have the same $\omega$-vertex say x and $v_p$ have the $\omega$-vertex say y. Then < $v_m, v_n, x, v_p, y$> = $H_1$, a contradiction.

**Case 3.2:** All three have the same $\omega$-vertex say x and the same $\omega+1$-vertex say y.
Now xy $\in$ E(G) (else x has another $\omega+1$-vertex say z and <$v_t, v_{t+1}, y, x, z$> = $H_1$). W.l.g. let $xv_i \notin$ E(G). Then $v_i$ has no $\omega$-vertex (else if z is the $\omega$-vertex of $v_i$, then <$v_t, v_{t+1}, x, v_i, z$> = $H_1$). Thus i > k. Again each of $v_t, v_{t+1}, v_k$ has another i-vertex (else color $v_i$ by $\omega$, $v_j$ by i and $v_l$ by l+1 for l < j where j $\in$ {t, t+1, k}). Let z be an i-vertex s.t. z $\neq v_i$ and $zv_t \in$ E(G). As G is $K_{1,3}$-free, xz $\in$ E(G) and x has at the most two i-vertices. Hence at least two of $v_t, v_{t+1}, v_k$ have the same i-vertex other than $v_i$. If {m, n, p} = {t, t+1, k} and say $v_m$ has a different i-vertex say w, then <$v_n, v_p, x, v_i, w$> = $H_1$>, a contradiction. Thus all of $v_t, v_{t+1}, v_k$ have the same i-vertex say z and xz $\in$ E(G). As <Q> is a maximum clique, clearly $\exists$ i, j s.t. $xv_i, xv_j \notin$ E(G) and zw $\notin$ E(G) where w is the other j-vertex of $v_t, v_{t+1}, v_k$. Clearly $zv_j$, $wv_i \in$ E(G). Now if $yv_i \notin$ E(G), then as before yz $\in$ E(G). But then if $yv_j \notin$ E(G), yz $\in$ E(G) and <x, y, w, z, $v_i$> = $H_1$ and if $yv_j \in$ E(G), then <y, z, x, $v_j, v_i$> = $H_1$, a contradiction. Hence $yv_j, yv_i \in$ E(G). Clearly $\exists$ s s.t. $yv_s \notin$ E(G). As before $xv_s \in$ E(G). But then <$v_i, v_j, y, v_s, x$> = $H_2$, a contradiction.

This proves the theorem.

**Examples to show that the Upper Bound is Tight:**
Let $G = C_{2n+1}$, $n > 1$. Then G is $\{K_{1,3}, H_1, H_2\}$-free, $\omega(G) = 2$, and $\chi(G) = 3 = \omega(G)+1$.

**Examples to show that $K_{1,3}$, $H_1$, $H_2$ are Necessary:**
- Mycielski graphs with $\chi(G) = k \geq 4$, have only $K_{1,3}$-induced, are triangle-free, hence have $\omega(G) = 2 < \chi(G) - 1$.
- Let $H = C_5$. Construct G from H by replacing each vertex by $K_m$. Then G is $\{K_{1,3}, H_2\}$-free, G has $H_1$ induced, $\omega(G) = 2m$, and $\chi(G) = \omega(G)+\lceil \frac{m}{2} \rceil > \omega(G)+1$ for $m \geq 3$.
- Let $G = \sum_{1}^{m} C_5$. Then G is $\{K_{1,3}, H_1\}$-free, G has $H_2$ induced, $\omega(G) = 2m$, and $\chi(G) = 3m > \omega(G)+1$ for $m > 1$.